\documentclass[10pt,twoside]{article}
\usepackage{amsfonts}
\usepackage{fancyhdr}
\usepackage{titlesec}
\usepackage{cite}
\usepackage{ifthen}
\usepackage{amssymb}
\usepackage{fancyhdr}
\usepackage{titlesec}
\usepackage[arrow,matrix]{xy}
\usepackage{pifont}
\usepackage{stmaryrd}
\usepackage{setspace}
\usepackage{indentfirst}
\usepackage{amsmath,amssymb,amscd,bbm,amsthm,mathrsfs,dsfont}
\usepackage{graphicx}
\input amssym.def

\def\setzero{\setcounter{equation}{0}}

\titleformat{\section}{\centering\large\bfseries}{\S\arabic{section}}{1em}{}
\newboolean{first}
\setboolean{first}{true}

\begin{document}

\setlength\abovedisplayskip{2pt}
\setlength\abovedisplayshortskip{0pt}
\setlength\belowdisplayskip{2pt}
\setlength\belowdisplayshortskip{0pt}

\title{\bf \Large On existence and uniqueness of solutions to
uncertain backward stochastic differential equations\author{FEI Wei-yin}\date{}} \maketitle
 \footnote{}
 \footnote{MR Subject Classification: 60H10, 94D05.}
%%%%%%%%%%%%%%%%%关键词至少三个及以上
 
 \footnote{ Supported by National Natural Science Foundation of China\,(71171003, 71210107026),
Anhui Natural Science Foundation\,(10040606003), and Anhui Natural
Science Foundation of Universities\,(KJ2012B019, KJ2013B023).}
\begin{center}
\begin{minipage}{135mm}

{\bf \small Abstract}.\hskip 2mm {\small This paper is concerned with a class of uncertain backward stochastic differential equations (UBSDEs) driven by both an $m$-dimensional Brownian motion and a $d$-dimensional canonical process with uniform Lipschitzian coefficients. Such equations can be useful in modelling hybrid systems, where the phenomena are simultaneously subjected
to two kinds of uncertainties: randomness and uncertainty. The solutions of UBSDEs are
the uncertain stochastic processes. Thus, the existence and uniqueness of solutions to UBSDEs with Lipschitzian coefficients are proved.}
\end{minipage}
\end{center}

\section{Introduction} \setzero

%\vskip12pt

Randomness is a basic type of objective uncertainty, and probability theory is a branch of
mathematics for studying the behavior of random phenomena. The study of probability theory was started by
Pascal and Fermat in 1654, and an axiomatic foundation of probability theory given by Kolmogorov in 1933.
The concept of fuzzy set was initiated by Zadeh \cite{Za1} via membership function in 1965.
In order to measure a fuzzy event, Zadeh \cite{Za2} introduced the theory of possibility. Moreover, fuzzy random variables are mathematical descriptions for fuzzy stochastic phenomena (i.e., a mixture of fuzziness
and randomness) and can be defined in several ways on the
basis of probability theory and fuzzy mathematics. The concept of fuzzy random variables are introduced by Kwakernaak \cite{Kw1, Kw2} and Puri and Ralescu \cite{PR}. Furthermore, the theory of fuzzy-valued (or set-valued) random systems was also investigated by many researchers, such as Fei \cite{FeiINS05,FeiINS,FLF,FLZ, FeiU, FeiNA, FX}, Li and Guan \cite{LG}, Malinowski et al. \cite{MalINS, MMS} and references therein.

However, some information and knowledge are usually represented by human language like ``about 100km'',
``roughly 80kg'', ``low speed'', ``middle age'', and ``big size''. A lot of surveys show
that in the real life imprecise quantities behave neither like randomness nor like fuzziness. When the sample size is too small (even no-sample) to estimate a probability
distribution, we have to invite some domain experts to evaluate their belief
degree on which each event will occur. Since human beings usually overestimate
unlikely events, the belief degree may have much larger variance than the
real frequency. Perhaps some people think that the belief degree is subjective
probability. However, it is inappropriate because probability theory may
lead to counterintuitive results in this case. In order to distinguish it from randomness, we name this phenomenon uncertainty. How do we understand
uncertainty? How do we model uncertainty? In order to answer those questions,
an uncertainty theory is founded in 2007 by Liu \cite{Liu2012}, which then becomes a branch
of mathematics for modeling human uncertainty.

We know that the additivity axiom of classical measure theory has
been challenged by many mathematicians. The earliest challenge was
from the theory of capacities by Choquet \cite{Ch} in which
monotonicity and continuity axioms are assumed.  For this reason,
an uncertainty theory founded in Liu \cite{Liu2012} gives  a new system of axioms based on normality,
self-duality, countable subadditivity and product measure.

Differential equations have been widely applied in physics, engineering, biology, economics and other fields.
With the development of science and technology, practical problems require more and more accurate description.
A wide range of uncertainties are added to the differential equation system, thus produce stochastic differential
equations, fuzzy differential equations and fuzzy stochastic differential equations.
Furthermore, uncertain differential equation, a type of differential equations driven by canonical process,
was defined by Liu \cite{Liu2012} in 2007. Chen and Liu \cite{Chen} present an existence and uniqueness theorem of solution for
uncertain differential equation under Lipschitz condition and linear growth condition. Zhu \cite{Zhu} investigates the uncertain optimal control with application to a portfolio selection model. In  Ge and Zhu \cite{GZ}, a necessary condition of optimality for uncertain optimal control problem is provided, where the existence and uniqueness of solutions to a backward uncertain differential equation is proven. The neutral uncertain delay differential equations and almost sure stability for uncertain differential equations are discussed in Liu et al. \cite{LF, LFS}. The optimal control of uncertain stochastic systems with
Markovian switching and its applications to portfolio decisions is investigated in \cite{FeiO}.

In the investigations of stochastic dynamic systems, the linear backward stochastic differential equations for the adjoint process in optimal stochastic control were early explored in Kushner \cite{Ku}, Yong and Zhou \cite{YZ} and references therein. In 1990, the adapted solution of a backward stochastic differential equation was successfully solved in Pardoux and Peng \cite{PP}, which started a new field of study involved in both the mathematical theory and the applications such as control theory, biology, engineering, economics and finance etc. (see e.g. \cite{CE,EPQ,FP,FSM,Fretire,FXZ}). Later, the adapted solution of backward stochastic differential equations with non-Lipschitzian coefficients was studied in Mao \cite{Mao}.

After further investigation, however, we find that, for many real important problems, one can assume that the system under consideration includes both the randomness and the Liu's uncertainty. In fact, the notion and its properties of uncertain random variables are studied by Liu \cite{LiuY}. Hence it might be appropriate that  uncertain stochastic dynamic systems are characterized by UBSDEs disturbed by both a Wiener process and a canonical process. In this paper, we first formulate this UBSDEs. Then the existence and uniqueness of solutions to UBSDEs is proven by classical martingale presentation theorem in stochastic calculus and a Picard type iteration. Our results will be meaningful to developing  the theory and applications of UBSDEs further.

The rest of the paper is organized as follows.
Some preliminary concepts of uncertainty theory are recalled in Section 2. And the notion of the UBSDE is also formulated.
Several key propositions are proven in Section 3. In Section 4, the proof of the main theorem is completed. Finally, conclusions are made in Section 5.

\section{Preliminaries}\setzero

Throughout this paper,  let $(\Omega,{\cal F},\{{\cal F}_t\}_{t\in [0,T]}, P)$ be a complete, filtered probability space where the sub-$\sigma$-field family $({\cal F}_t, t\in [0,T])$ of $\cal F$
satisfies the usual conditions and $T\in[0,\infty)$ is a time horizon. The Brownian filtration $\{{\cal F}_t\}_{t\in[0,T]}$ is generalized by $\sigma(B_s:s\leq t)$ and $P$-null sets of $\cal F$, ${\cal F}_T={\cal F}$, where an $m$-dimensional Wiener process $B_t=(B_t^1,\cdots B_t^m)^\top$ is defined on the probability space $(\Omega, {\cal F}, P)$. The related properties on probability space refer to \cite{KS,Ok}.

Let $(\Gamma, {\cal L}, {\cal M})$ be an uncertainty space described in Liu \cite{Liu2012} where normality,
self-duality, countable subadditivity and product measure axioms are fulfilled.
Now we define a filtration which is the sub-$\sigma$-field family $({\cal L}_t, t\in [0,T])$ of $\cal L$
satisfying the usual conditions. The canonical process filtration $\{{\cal L}_t\}_{t\in[0,T]}$ is generalized by $\sigma(C_s:s\leq t)$ and $\cal M$-null sets of $\cal L$, ${\cal L}_T={\cal L}$, where a $d$-dimensional canonical process $C_t=(C_t^1,\cdots C_t^d)^\top$ is defined on the uncertainty space $(\Gamma, {\cal L}, {\cal M})$.

In order to discuss the uncertain stochastic systems, we need to construct a filtered uncertain probability space $(\Gamma\times\Omega, {\cal L}\otimes{\cal F}, ({\cal L}_t\otimes{\cal F}_t)_{t\in[0,T]},{\cal M}\times{P})$ on which we can define the related concepts as follows.

\vskip12pt

{\noindent\bf  Definition 2.1.} (i) An uncertain random variable is a measurable function $\xi\in{\mathbb R}^p$ (resp. ${\mathbb R}^{p\times m}$) from an uncertainty probability space $(\Gamma\times\Omega, {\cal L}\otimes{\cal F}, {\cal M}\times{P})$ to the set in ${\mathbb R}^p$ (resp. ${\mathbb R}^{p\times m}$), i.e., for any Borel set $A$ in ${\mathbb R}^p$ (resp. ${\mathbb R}^{p\times m}$), the set $\{\xi\in A\}=\{(\gamma,\omega)\in \Gamma\times\Omega: \xi(\gamma,\omega)\in A\}\in {\cal L}\otimes{\cal F}$.

(ii) The expected value of an uncertain random variable $\xi$ is defined by $${\mathbb E}[\xi]={\mathbb E}_{P}[{\mathbb E}_{\cal M}[\xi]]\stackrel{\triangle}{=}\int_\Omega\left[\int_0^{+\infty}{\cal M}\{\xi\geq r\}dr\right]P(d\omega)-\int_\Omega\left[\int_{-\infty}^0{\cal M}\{\xi\leq r\}dr\right]P(d\omega),$$
where ${\mathbb E}_{P}$ and ${\mathbb E}_{\cal M}$ denote the expected values under the uncertainty space and the probability space, respectively.

 \vskip12pt

 Obviously, if both $a$ and $b$ are constant, then ${\mathbb E}[aC_t+bB_t]=0$, where $C_t$ and $B_t$ are a scalar canonical process and a  one-dimensional Wiener process (Brownian motion), respectively. Notice that our definition on uncertain random variables is slightly different from the one in Liu \cite{LiuY}, where an uncertain random variable is roughly a function from a probability space to the set of uncertain variables.
 \vskip12pt

{\noindent \bf Definition 2.2.}  An uncertain process $X(t)\in{\mathbb R}^p$ (resp. ${\mathbb R}^{p\times m}$) is a measurable function from $[0,T]\times\Gamma$ to the set in ${\mathbb R}^p$ (resp. ${\mathbb R}^{p\times m}$), i.e., for each $t\in [0,T]$ and any Borel set $A$
in ${\mathbb R}^p$ (resp. ${\mathbb R}^{p\times m}$), the set $\{X(t)\in A\}=\{\gamma\in \Gamma| X(t,\gamma)\in A\}\in {\cal L}$.

\vskip12pt

The concepts and properties of the canonical process and other uncertain processes refer to Chapter 12 Liu \cite{Liu2012}. For $x\in {\mathbb R}^p, |x|$ denote its Euclidean norm. An element $y\in {\mathbb R}^{p\times m}$ will be considered as a $p\times m$ matrix; note that its Euclidean norm is given by $|y|=\sqrt{Tr(yy^\top)},$ and $(y,z)=Tr(yz^\top), z\in {\mathbb R}^{p\times m}$. In what follows, we give the notation of uncertain stochastic processes.

 \vskip12pt

{\noindent \bf Definition 2.3.}  (i) A hybrid process $X(t)$ is called an uncertain stochastic process if for each $t\in [0,T]$, $X(t)$ is an uncertain random variable. An uncertain stochastic process $X(t)$ is called continuous if the sample paths of $X(t)$ are all continuous
functions of $t$ for almost all $(\gamma, \omega)\in \Gamma\times\Omega$.

(ii) An uncertain stochastic process $X(t)$ is called ${\cal F}_t$-adapted if $X(t,\gamma)$ is ${\cal F}_t$-measurable for all $t\in[0,T], \gamma\in \Gamma$. Moreover, an uncertain stochastic process $X(t)$ is called ${\cal L}_t\otimes{\cal F}_t$-adapted (or adapted) if $X(t)$ is ${\cal L}_t\otimes{\cal F}_t$-measurable for all $t\in[0,T]$.

(iii) An uncertain stochastic process is called progressively measurable if it is measurable with respect to the $\sigma$-algebra
$$\Im({\cal L}_t\otimes{\cal F}_t))=\{A\in{\cal B}([0,T])\otimes{\cal L}\otimes{\cal F}: A\cap([0,t]\times{\Gamma\times\Omega})\in {\cal B}([0,t])\otimes{\cal L}_t\otimes{\cal F}_t\}.$$
 Moreover, an uncertain stochastic process $X(t):\Gamma\times\Omega\rightarrow {\mathbb R}^p$ (resp. $X(t):\Gamma\times\Omega\rightarrow {\mathbb R}^{p\times m}$) is called $L^2$-progressively measurable if it is progressively measurable and satisfies $E\left[\int_0^T|X(t)|^2dt\right]<\infty$. ${\mathbb M}^2(0,T; {\mathbb R}^p) $ (resp. ${\mathbb M}^2(0,T;{\mathbb R}^{p\times m })$ denote the set of $L^2$-progressively measurable uncertain random processes.

 \vskip12pt
{\noindent \bf Definition 2.4.} (It\^o-Liu integral) Let $X(t)=(Y(t),Z(t))^\top$ be an uncertain stochastic process, where $Y(t)\in{\mathbb R}^{p\times m}$ and  $Z(t)\in{\mathbb R}^{p\times d}$.  For any partition of
closed interval $[a, b]$ with $a = t_1 < t_2 < \cdots < t_{N+1} = b,$ the mesh is written as $\Delta=\max\limits_{1\leq i\leq N}|t_{i+1}-t_i|.$
Then the It\^o-Liu integral of $X(t)$ with respect to $(B_t, C_t)$ is defined as follows,
$$\int_a^bX(s)^\top d
\left(\begin{array}{ll}
&B_s\\
&C_s
\end{array}\right)
=\lim\limits_{\Delta\rightarrow0}\sum\limits_{i=1}^N(Y({t_i})(B_{t_{i+1}}-B_{t_i})+Z({t_i})(C_{t_{i+1}}-C_{t_i}))$$
provided that it exists in mean square and is an uncertain random variable, where $C_t$ and $B_t$ are a $d$-dimensional canonical process and an $m$-dimensional Wiener process, respectively.
 In this case, $X(t)$ is
called It\^o-Liu integrable. Specially, when $Y(t)\equiv0$, $X(t)$ is called Liu integrable.

\vskip12pt

The following It\^o-Liu formula for the case of multi-dimensional uncertain stochastic processes (see e.g. Fei \cite{FeiO}) is given.

\vskip12pt

{\noindent \bf Theorem 2.5.} {\sl Let
$$
\begin{array}{ll}
B&=(B_t)_{0\leq t\leq T}=(B^1_t, \cdots,
B_t^m)^\top_{0\leq t\leq T}\mbox{ and}\\
 C&=(C_t)_{0\leq t\leq T}=(C^1_t, \cdots, C_t^d)^\top_{0\leq t\leq T}
\end{array}
$$
 be an $m$-dimensional
standard Wiener process and a $d$-dimensional canonical process, respectively.
Assume that uncertain stochastic processes $X_1(t), X_{2}(t), \cdots, X_{p}(t)$ are given by
$${d} X_{k}(t)=u_{k}(t){d}t+\sum\limits_{l=1}^mv_{kl}(t){d}B_{t}^l+\sum\limits_{l=1}^dw_{kl}(t){d}C_{t}^l,\quad k=1, \cdots, p, $$
where $u_{k}(t)$ are all absolute integrable uncertain stochastic processes, $v_{kl}(t)$ are all square integrable
uncertain stochastic processes and $w_{kl}(t)$ are all Liu integrable uncertain stochastic processes. For $k,l=1, \cdots, p$,  let $\frac{\partial G}{\partial t}(t, x_1, \cdots, x_p)$, $\frac{\partial G}{\partial x_k}(t, x_1, \cdots, x_p)$ and $\frac{\partial^2 G}{\partial x_kx_l}(t, x_1, \cdots, x_p)$ be continuously functions. Then we have
$$
\begin{array}{ll}
&d G(t, X_{1}(t), \cdots, X_{p}(t))\\
&=\frac{\partial G}{\partial t}(t, X_{1}(t), \cdots, X_{p}(t)){d }t+\sum\limits_{k=1}^p\frac{\partial G}{\partial x_k}(t, X_{1}(t), \cdots, X_{p}(t)){d}X_{k}(t)\\
&\quad+\frac{1}{2}\sum\limits_{k=1}^p\sum\limits_{l=1}^p\frac{\partial^2 G}{\partial x_k\partial x_l}(t, X_{1}(t), \cdots, X_{p}(t)){d}X_{k}(t){d}X_{l}(t),
\end{array}
$$
where ${d}B_{t}^k{d}B_{t}^l=\delta_{kl}{d}t, {d}B_{t}^k{d}t={d}C_{t}^\imath{d}C_{t}^\jmath={d}C_{t}^\imath{d}t={d}B_{t}^k{d}C_{t}^\imath=0,$ for
 $k,l=1,\cdots, m, \imath,\jmath=1, \cdots, d.$ Here
$$
\begin{array}{ll}
\delta_{kl}=\left\{\begin{array}{ll}
0,& {\rm if}\quad k\neq l,\\
1, & \rm {otherwise}.
\end{array}
\right.
\end{array}
$$
}

\vskip12pt

In what follows, we consider
the following uncertain backward stochastic differential equation (UBSDE), for $t\in[0,T]$,
$$
\begin{array}{ll} dX(t)=&f(t, X(t), Y(t))dt+g(t,
X(t),Y(t))dC(t)\\
&+h(t,X(t),Y(t))dB(t), \quad
X(T)=\xi,
\end{array}
$$
or, in the integral form,
$$
\begin{array}{ll} X(t)+&\int_t^Tf(s, X(s), Y(s))ds+\int_t^Tg(s,
X(s),Y(s))dC(s)\\
&+\int_t^Th(s,X(s),Y(s))dB(s)=\xi,
\end{array}
 \eqno(2.1)
$$
where
$$
\begin{array}{ll}
f=&(f_1,\cdots,f_p)^\top: \Gamma\times\Omega\times[0,T]\times{\mathbb R}^p\times
{\mathbb
R}^{p\times m}\rightarrow {\mathbb R}^p\\
&\quad\mbox{being}~{\cal P}\otimes{\cal B}_p\otimes{\cal B}_{p\times m}/{\cal B}_{p}~\mbox{measurable},\\
 g=&(g_{kl})_{p\times d}: \Gamma\times\Omega\times[0,T]\times{\mathbb R}^p\times{\mathbb
R}^{p\times m}\rightarrow {\mathbb
R}^{p\times d}\\
&\quad\mbox{being}~{\cal P}\otimes{\cal B}_p\otimes{\cal B}_{p\times m}/{\cal B}_{p\times d}~\mbox{measurable},\\
h=&(h_{kl})_{p\times m}: \Gamma\times\Omega\times[0,T]\times{\mathbb R}^p\times{\mathbb
R}^{p\times m}\rightarrow {\mathbb
R}^{p\times m}\\
&\quad\mbox{being}~{\cal P}\otimes{\cal B}_p\otimes{\cal B}_{p\times m}/{\cal B}_{p\times m}~\mbox{measurable}.
\end{array}
$$
Here, $\cal P$ denotes the $\sigma$-algebra of progressively measurable subsets of $\Gamma\times\Omega\times[0,T]$.

Our main result will be an existence and uniqueness result for an adapted pair $\{X(t),Y(t); t\in[0,T]\}$ which solves (2.1). In order to guarantee the existence and uniqueness of solution to UBSDEs, we suppose that  $h(\cdot)$ satisfies a rather restrictive assumption which implies in particular that the mapping $y\rightarrow h(s, x, y)$ is a bijection for any $(\gamma,\omega,s, x)$.

In terms of definitions of the expectation operator ${\mathbb E}[\cdot]$ and It\^o-Liu uncertain stochastic integral, it is easy to see that, for $\forall a, b\in[0, T], X\in{\mathbb M}^2(0,T;{\mathbb R}^{p\times d})$, and $ Y\in{\mathbb M}^2(0,T;{\mathbb R}^{p\times m})$,
$${\mathbb E}\left[\int_a^bX(t)dC(t)\right]=0~\mbox{and }~{\mathbb E}\left[\int_aY(t)dB(t)\right]=0$$
which are frequently used in the next section.

\vskip12pt

\section{Results on simplified versions of UBSDE (2.1)}

For the aim of the study of UBSDE (2.1), in this section, we provide three simplified versions of equation (2.1).  First, we consider the following UBSDE

$$
 X(t)+\int_t^Tf(s)ds+\int_t^Tg(s)dC(s)+\int_t^T[h(s)+Y(s)]dB(s)=\xi,~0\leq t\leq T.\eqno(3.1)
$$

\vskip12pt

{\noindent\bf Proposition 3.1.} {\sl Given $\xi\in L^2(\Gamma\times\Omega, {\cal L}\otimes{\cal F},  {\cal M}\times P; {\mathbb R}^p)$,
$f\in {\mathbb M^2(0,T;{\mathbb R}^p)}$, $g\in{\mathbb M}^2(0,T;{\mathbb R}^{p\times d})$ and $h\in {\mathbb M}^2(0,T; {\mathbb R}^{p\times m})$, there exists a unique pair $(X, Y)\in {\mathbb M}^2(0,T;{\mathbb R}^p)\times {\mathbb M}^2(0,T;{\mathbb R}^{p\times m})$ such that UBSDE (3.1) holds.
}

\vskip12pt

{\noindent\it Proof.} For any fixed $\gamma\in {\Gamma}$, we define
$$X(t,\gamma)={\mathbb E}_{P}\left[\xi(\gamma)-\int_t^Tf(s,\gamma)ds-\int_t^Tg(s,\gamma)dC(s)|{\cal F}_t\right],~0\leq t\leq T.$$
Then
$$\tilde{X}(t,\gamma)={\mathbb E}_{P}\left[\xi(\gamma)-\int_0^Tf(s,\gamma)ds-\int_0^Tg(s,\gamma)dC(s)|{\cal F}_t\right],~0\leq t\leq T$$
 is an ${\cal F}_t$-adapted martingale on the filtered probability space $(\Omega,{\cal F},\{{\cal F}_t\}_{t\in [0,T]}, P)$ with $\tilde{X}(0,\gamma)=X(0,\gamma)$ for almost all $\gamma\in \Gamma$.
By a well-known martingale representation theorem (see e.g. Karatzas and Shreve \cite{KS}, p.182 or $\O$ksendal \cite{Ok}, p.53), it follows from the assumptions that there exists $\bar{Y}\in {\mathbb M}(0,T, {\mathbb R}^{p\times m})$ such that
$${\mathbb E}\left[\xi(\gamma)-\int_0^Tf(s,\gamma)ds-\int_0^Tg(s,\gamma)dC(s)|{\cal F}_t\right]=X( 0,\gamma)+\int_0^t\bar{Y}(s,\gamma)dB(s),$$
where the operator ${\mathbb E}[\cdot|{\cal F}_t]$ denotes the conditional expectation of a stochastic process with respect to the filtration ${\cal F}_t$. We now define
$Y(t)=\bar{Y}(t)-h(t), 0\leq t\leq T$. It is easily seen that the constructed pair $(X,Y)$ being adapted solves UBSDE (3.1).

For the proof of uniqueness, let $(X_1, Y_1)$ and $(X_2,Y_2)$ be two solutions to UBSDE (3.1). From It\^o-Liu formula (see Theorem 2.5) applied to $|X_1(s)-X_2(s)|^2$ from $s=t$ to $T$ we have
$$
\begin{array}{ll}
&|X_1(t)-X_2(t)|^2+\int_t^T|Y_1(s)-Y_2(s)|^2ds=-2\int_t^T(X_1(s)-X_2(s),[Y_1(s)-Y_2(s)]dB(s)).
\end{array}
$$
Hence we deduce
$$
{\mathbb E}|X_1(t)-X_2(t)|^2+\int_t^T{\mathbb E}|Y_1(s)-Y_2(s)|^2ds=0,
$$
which shows $X_1(t)=X_2(t), Y_1(t)=Y_2(t), {\cal M}\times P$-a.e. Thus the proof is complete. $\Box$

\vskip12pt

We now consider the UBSDE

$$
 X(t)+\int_t^Tf(s, Y(s))ds+\int_t^Tg(s,Y(s))dC(s)+\int_t^T[h(s)+Y(s)]dB(s)=\xi,~0\leq t\leq T,\eqno(3.2)
$$
where
$f: \Gamma\times\Omega\times[0,T]\times
{\mathbb
R}^{p\times m}\rightarrow {\mathbb R}^p$ is  ${\cal P}\otimes{\cal B}_{p\times m}/{\cal B}_{p}$ measurable and $g: \Gamma\times\Omega\times[0,T]\times
{\mathbb
R}^{p\times m}\rightarrow {\mathbb R}^{p\times d}$ is  ${\cal P}\otimes{\cal B}_{p\times m}/{\cal B}_{p\times d}$ measurable with the property that
$$f(\cdot, 0)\in {\mathbb M}^2(0,T;{\mathbb R}^p), \quad g(\cdot, 0)\in {\mathbb M}^2([0,T]; {\mathbb R}^{p\times d})\eqno(3.3)$$
and there exists $c>0$ such that
$$|f(t, y_1)-f(t,y_2)|\vee |g(t, y_1)-g(t,y_2)|\leq c|y_1-y_2|, \eqno(3.4)$$
for any $y_1,y_2\in{\mathbb R}^{p\times m}$, and $(\gamma, \omega,t)$ a.e. Note that (3.3) and (3.4) imply that $f(\cdot, Y(\cdot))\in {\mathbb M}^2(0,T;{\mathbb R}^p)$ and $g(\cdot, Y(\cdot))\in {\mathbb M}^2(0,T; {\mathbb R}^{p\times d})$ whenever $Y\in {\mathbb M}^2(0,T;{\mathbb R}^{p\times m}).$

\vskip12pt

{\noindent\bf Proposition 3.2.} {\sl In  UBSDE (3.2), let $\xi\in L^2(\Gamma\times\Omega, {\cal L}\otimes{\cal F},  {\cal M}\times P; {\mathbb R}^p)$,
$g\in{\mathbb M}^2(0,T;{\mathbb R}^{p\times d})$, $h\in {\mathbb M}^2(0,T;{\mathbb R}^{p\times m})$. If $f: \Gamma\times\Omega\times[0,T]\times{\mathbb R}^{p\times m}\rightarrow {\mathbb R}^p$ is a mapping satisfying the above requirements, in particular (3.4), then there exists a unique pair $(X, Y)\in {\mathbb M}^2(0,T;{\mathbb R}^p)\times {\mathbb M}^2(0,T;{\mathbb R}^{p\times m})$ such that
 UBSDE (3.1) holds.
}

\vskip12pt

{\noindent\it Proof.} {\underline{Uniqueness}}.  Let $(X_1, Y_1)$ and $(X_2,Y_2)$ be two solution to UBSDE (3.2). From It\^o-Liu formula (see Theorem 2.5) applied to $|X_1(s)-X_2(s)|^2$ from $s=t$ to $T$ we have
$$
\begin{array}{ll}
&|X_1(t)-X_2(t)|^2+\int_t^T|Y_1(s)-Y_2(s)|^2ds\\
&=-2\int_t^T(f(s,Y_1(s))-f(s, Y_2(s)),X_1(s)-X_2(s))ds\\
&\quad-2\int_t^T(X_1(s)-X_2(s), [g(s,Y_1(s))-g(s, Y_2(s))]dC(s))\\
&\quad-2\int_t^T(X_1(s)-X_2(s),[Y_1(s)-Y_2(s)]dB(s)),
\end{array}
$$
from which and (3.4) we deduce
$$
\begin{array}{ll}
&{\mathbb E}|X_1(t)-X_2(t)|^2+{\mathbb E}\int_t^T|Y_1(s)-Y_2(s)|^2ds\\
&=-2{\mathbb E}\int_t^T(f(s,Y_1(s))-f_2(s, Y_2(s)),X_1(s)-X_2(s))ds\\
&\leq \frac{1}{2}{\mathbb E}\int_t^T|Y_1(s)-Y_2(s)|^2ds+2c^2\int_t^TE|X_1(s)-X_2(s)|^2ds,
\end{array}
$$
which shows, from Gronwall inequality, $X_1(t)=X_2(t), Y_1(t)=Y_2(t)$ for almost all $\gamma\times\omega\in \Gamma\times\Omega$.

{\underline{Existence}}. Due to Proposition 3.1, we define an approximating sequence by a Picard iteration. Let $Y_0(t)\equiv0$, and $\{(X_n(t),Y_n(t)); 0\leq t\leq T\}_{n\geq1}$ be a sequence in ${\mathbb M}^2(0, T;{\mathbb R}^p)\times{\mathbb M}^2(0, T;{\mathbb R}^{p\times m})$ define recursively by
$$
 X_n(t)+\int_t^Tf(s, Y_{n-1}(s))ds+\int_t^Tg(s,Y_{n-1}(s))dC(s)+\int_t^T[h(s)+Y_n(s)]dB(s)=\xi.\eqno(3.5)
$$
By Theorem 2.5 and same inequalities as above, we have $(K=2c^2)$
$$
\begin{array}{ll}
&{\mathbb E}|X_{n+1}(t)-X_{n}(t)|^2+{\mathbb E}\int_t^T|Y_{n+1}(s)-Y_n(s)|^2ds\\
&\leq \frac{1}{2}{\mathbb E}\int_t^T|Y_n(s)-Y_{n-1}(s)|^2ds+K\int_t^T{\mathbb E}|X_{n+1}(s)-X_n(s)|^2ds.
\end{array}
\eqno(3.6)
$$
Set $\varphi_n(t)=\int_t^T{\mathbb E}|X_n(s)-X_{n-1}(s)|^2ds$ and $\psi_n(t)={\mathbb E}\int_t^T|Y_n(s)-Y_{n-1}(s)|^2ds$, for $n\geq 1$ $(X_0(t)\equiv0$).
From (3.6), we have
$$-\frac{d}{dt}\left(\varphi_{n+1}(t)e^{Kt}\right)+e^{Kt}\psi_{n+1}(t)\leq \frac{1}{2}e^{Kt}\psi_n(t).\eqno(3.7)$$
Integrating from $t$ to $T$, we have
 $$\varphi_{n+1}(t)+\int_t^Te^{K(s-t)}\psi_{n+1}(s)ds\leq \frac{1}{2}\int_t^Te^{K(s-t)}\psi_n(s)ds,$$
which implies that
$$\int_0^Te^{Kt}\psi_{n+1}(t)dt\leq 2^{-n}\tilde{c}e^{KT}$$
with $\tilde{c}={\mathbb E}\int_0^T|Y_1(t)|^2dt=\sup_{0\leq t\leq T}\psi_1(t)$; but also then
$$\psi_{n+1}(0)\leq 2^{-n}\tilde{c}e^{KT}.\eqno(3.8)$$
On the other hand, from (3.7) and $(d/dt)\varphi_{n+1}(t)\leq0$ we get
$$\psi_{n+1}(0)\leq K\varphi_{n+1}(0)+\frac{1}{2}\psi_n(0)\leq 2^{-n}\tilde{K}+\frac{1}{2}\psi_n(0),$$
where $\tilde{K}=\tilde{c}Ke^{KT}.$ Thus by iterating we get
$$\psi_{n+1}(0)\leq 2^{-n}(n\tilde{K}+\psi_1(0)).\eqno(3.9)$$
Since the square roots of the right hand sides of (3.8) and (3.9) are summable series, we know that $\{X_n\}$ (resp. $\{Y_n\}$) is a Cauchy sequence in ${\mathbb M}^2(0,T;{\mathbb R}^p)$ (resp. ${\mathbb M}^2(0,T;{\mathbb R}^{p\times m}$). Hence from (3.5), $X_n$ is also a Cauchy sequence in $L^2(\Gamma\times\Omega; C([0,T];{\mathbb R}^p))$, and passing to the limit in (3.5) as $n\rightarrow\infty$, we obtain that the pair $(X,Y)$ defined by
$$X=\lim\limits_{n\rightarrow\infty}X_n,\quad Y=\lim\limits_{n\rightarrow\infty}Y_n$$
solves the UBSDE (3.2). Thus the proof is complete. \hfill$\Box$

\vskip12pt
 We now study the UBSDE

$$
 \begin{array}{ll}
 &X(t)+\int_t^Tf(s,X(s), Y(s))ds+\int_t^Tg(s,X(s),Y(s))dC(s)\\
 &\quad+\int_t^T[h(s, X(s))+Y(s)]dB(s)=\xi,~0\leq t\leq T,
 \end{array}
 \eqno(3.10)
 $$
where
$$
\begin{array}{ll}
f=&(f_1,\cdots,f_p)^\top: \Gamma\times\Omega\times[0,T]\times{\mathbb R}^p\times
{\mathbb
R}^{p\times m}\rightarrow {\mathbb R}^p\\
&\quad\mbox{being}~{\cal P}\otimes{\cal B}_p\otimes{\cal B}_{p\times m}/{\cal B}_{p}~\mbox{measurable},\\
 g=&(g_{kl})_{p\times d}: \Gamma\times\Omega\times[0,T]\times{\mathbb R}^p\times{\mathbb
R}^{p\times m}\rightarrow {\mathbb
R}^{p\times d}\\
&\quad\mbox{being}~{\cal P}\otimes{\cal B}_p\otimes{\cal B}_{p\times m}/{\cal B}_{p\times d}~\mbox{measurable},\\
h=&(h_{kl})_{p\times m}: \Gamma\times\Omega\times[0,T]\times{\mathbb R}^p\rightarrow {\mathbb
R}^{p\times m}\\
&\quad\mbox{being}~{\cal P}\otimes{\cal B}_p/{\cal B}_{p\times m}~\mbox{measurable}.
\end{array}
$$
Moreover, the following properties hold
$$f(\cdot, 0)\in {\mathbb M}^2(0,T;{\mathbb R}^p), ~ g(\cdot, 0)\in {\mathbb M}^2(0,T; {\mathbb R}^{p\times d})~\mbox{and}~ h(\cdot, 0)\in {\mathbb M}^2(0,T; {\mathbb R}^{p\times m}),$$
and there exists $c>0$ such that
$$
\begin{array}{ll}
&|f(t, x_1, y_1)-f(t,x_2, y_2)|\vee |g(t,x_1, y_1)-g(t,x_2,y_2)|\leq c(|x_1-x_2|+|y_1-y_2|),\\
& |h(t, x_1)-h(t,x_2)|\leq c|x_1-x_2|
\end{array}
\eqno(3.11)
$$
for all $x_1,x_2\in{\mathbb R}^p, y_1,y_2\in{\mathbb R}^{p\times m}, (\gamma,\omega,t)$-a.e. Note that (3.11) implies that $f(\cdot,X(\cdot), Y(\cdot))\in {\mathbb M}^2(0,T;{\mathbb R}^p)$, $g(\cdot, X(\cdot), Y(\cdot))\in {\mathbb M}^2(0,T; {\mathbb R}^{p\times d})$ and $h(\cdot, X(\cdot))\in {\mathbb M}^2(0,T; {\mathbb R}^{p\times m})$ whenever $X\in {\mathbb M}^2(0,T;{\mathbb R}^p)$ and $Y\in {\mathbb M}^2(0,T;{\mathbb R}^{p\times m}).$

\vskip12pt

{\noindent\bf Proposition 3.3.} {\sl In the UBSDE (3.10), let $\xi\in L^2(\Gamma\times\Omega, {\cal L}\otimes{\cal F},  {\cal M}\times P; {\mathbb R}^p)$,
$g\in{\mathbb M}^2(0,T;{\mathbb R}^{p\times d})$, $h\in {\mathbb M}^2(0,T;{\mathbb R}^{p\times m})$. If $f: \Gamma\times\Omega\times[0,T]\times{\mathbb R}^{p}\times {\mathbb R}^{p\times m}\rightarrow {\mathbb R}^p$ is a mapping satisfying the above requirements, in particular (3.11), then there exists a unique pair $(X, Y)\in {\mathbb M}^2(0,T;{\mathbb R}^p)\times {\mathbb M}^2(0,T;{\mathbb R}^{p\times m})$ such that UBSDE (3.10) holds.
}

\vskip12pt

{\noindent\it Proof.} {\underline{Uniqueness}}.  Let $(X_1, Y_1)$ and $(X_2,Y_2)$ be two solution to UBSDE (3.10). By similar argument as in the one in Proposition 3.2, we have
$$
\begin{array}{ll}
&|X_1(t)-X_2(t)|^2+\int_t^T|Y_1(s)-Y_2(s)|^2ds\\
&=-2\int_t^T(f(s,X_1(s),Y_1(s))-f(s, X_2(s),Y_2(s)),X_1(s)-X_2(s))ds\\
&\quad-2\int_t^T(X_1(s)-X_2(s), [g(s,X_1(s),Y_1(s))-g(s, X_2(s),Y_2(s))]dC(s))\\
&\quad-2\int_t^T(X_1(s)-X_2(s),[h(s,X_1(s))-h(s, X_2(s))+Y_1(s)-Y_2(s)]dB(s))\\
&\quad-2\int_t^T(h(s,X_1(s))-h(s,X_2(s)),[Y_1(s)-Y_2(s)])ds\\
&\quad-\int_t^T|h(s,X_1(s))-h(s, X_2(s))|^2ds.
\end{array}
$$
From (3.11) we deduce
$$
\begin{array}{ll}
&{\mathbb E}|X_1(t)-X_2(t)|^2+{\mathbb E}\int_t^T|Y_1(s)-Y_2(s)|^2ds\\
&=-2{\mathbb E}\int_t^T(f(s,X_1(s), Y_1(s))-f_2(s, X_2(s),Y_2(s)),X_1(s)-X_2(s))ds\\
&\quad-2{\mathbb E}\int_t^T(h(s,X_1(s))-h(s,X_2(s)),[Y_1(s)-Y_2(s)])ds\\
&\quad-{\mathbb E}\int_t^T|h(s,X_1(s))-h(s, X_2(s))|^2ds\\
&\leq \frac{1}{2}{\mathbb E}\int_t^T|Y_1(s)-Y_2(s)|^2ds+\bar{c}\int_t^TE|X_1(s)-X_2(s)|^2ds,
\end{array}
$$
where $\bar{c}$ is some constant. Thus from Gronwall inequality, the uniqueness easily is obtained.

{\underline{Existence}}. With the help of Proposition 3.2 we construct an approximating sequence by using a Picard iteration. Let $X_0(t)\equiv0$, and $\{(X_n(t),Y_n(t)); 0\leq t\leq T\}_{n\geq1}$ be a sequence in ${\mathbb M}^2(0, T;{\mathbb R}^p)\times{\mathbb M}^2(0, T;{\mathbb R}^{p\times m})$ defined recursively by
$$
\begin{array}{ll}
 X_n(t)&+\int_t^Tf(s,X_{n-1}(s), Y_n(s))ds+\int_t^Tg(s,X_{n-1}(s),Y_n(s))dC(s)\\
&+\int_t^T[h(s,X_{n-1})+Y_n(s)]dB(s)=\xi,~0\leq t\leq T.
\end{array}
\eqno(3.12)
$$
By Theorem 2.5 and same inequalities as above, we have
$$
\begin{array}{ll}
&|X_{n+1}(t)-X_n(t)|^2+\int_t^T|Y_{n+1}(s)-Y_n(s)|^2ds\\
&=-2\int_t^T(f(s,X_n(s),Y_{n+1}(s))-f(s, X_{n-1}(s),Y_n(s)),X_{n+1}(s)-X_n(s))ds\\
&\quad-2\int_t^T(X_{n+1}(s)-X_n(s), [g(s,X_{n}(s),Y_{n+1}(s))-g(s, X_{n-1}(s),Y_n(s))]dC(s))\\
&\quad-2\int_t^T(X_{n+1}(s)-X_n(s),[h(s,X_n(s))-h(s, X_{n-1}(s))+Y_{n+1}(s)-Y_n(s)]dB(s))\\
&\quad-2\int_t^T(h(s,X_n(s))-h(s,X_{n-1}(s)),[Y_{n+1}(s)-Y_n(s)])ds\\
&\quad-\int_t^T|h(s,X_n(s))-h(s, X_{n-1}(s))|^2ds,
\end{array}
$$
which, together with (3.11), shows
$$
\begin{array}{ll}
&{\mathbb E}|X_{n+1}(t)-X_{n}(t)|^2+{\mathbb E}\int_t^T|Y_{n+1}(s)-Y_n(s)|^2ds\\
&\leq \frac{1}{2}{\mathbb E}\int_t^T|Y_{n+1}(s)-Y_n(s)|^2ds\\
&\quad+c_1\left(\int_t^T{\mathbb E}|X_{n+1}(s)-X_n(s)|^2ds+\int_t^T{\mathbb E}|X_n(s)-X_{n-1}(s)|^2ds\right),
\end{array}
$$
where $c_1$ is a certain constant. Hence we have
$$
\begin{array}{ll}
&{\mathbb E}|X_{n+1}(t)-X_{n}(t)|^2+\frac{1}{2}{\mathbb E}\int_t^T|Y_{n+1}(s)-Y_n(s)|^2ds\\
&\leq c_1\left(\int_t^T{\mathbb E}|X_{n+1}(s)-X_n(s)|^2ds+\int_t^T{\mathbb E}|X_n(s)-X_{n-1}(s)|^2ds\right).
\end{array}
\eqno(3.13)
$$
Set $\varphi_n(t)=\int_t^T{\mathbb E}|X_n(s)-X_{n-1}(s)|^2ds$. It follows from (3.13) that
$$-\frac{d}{dt}\left(\varphi_{n+1}(t)\right)-c_1\varphi_{n+1}(t)\leq c_1\varphi_n(t),~\varphi_{n+1}(T)=0,$$
which implies that
$$\varphi_{n+1}(t)\leq c_1\int_t^Te^{c_1(s-t)}\varphi_n(s)ds.$$
Iterating that inequality, we have
$$\varphi_{n+1}(0)\leq \frac{(c_1e^{c_1T})^n}{n!}\varphi_1(0).$$
Thus together (3.13) we have that $\{X_n\}$ (resp. $\{Y_n\}$) is a Cauchy sequence in ${\mathbb M}^2(0,T;{\mathbb R}^p)$ (resp. ${\mathbb M}^2(0,T;{\mathbb R}^{p\times m}$). Hence from (3.12), $X_n$ is also a Cauchy sequence in $L^2(\Gamma\times\Omega; C([0,T];{\mathbb R}^p))$, and passing to the limit in (3.12) as $n\rightarrow\infty$, we obtain that the pair $(X,Y)$ defined by
$$X=\lim\limits_{n\rightarrow\infty}X_n,\quad Y=\lim\limits_{n\rightarrow\infty}Y_n$$
solves  UBSDE (3.10). Thus the proof is complete. \hfill$\Box$

\vskip12pt

\section{Existence and uniqueness of solutions to UBSDE (2.1)}

In this section, we consider UBSDE (2.1) as follows
$$
 \begin{array}{ll}
 X(t)&+\int_t^Tf(s,X(s), Y(s))ds+\int_t^Tg(s,X(s),Y(s))dC(s)\\
 &+\int_t^Th(s, X(s),Y(s))dB(s)=\xi,~0\leq t\leq T,
 \end{array}
 \eqno(4.1)
 $$
where
$$
\begin{array}{ll}
f=&(f_1,\cdots,f_p)^\top: \Gamma\times\Omega\times[0,T]\times{\mathbb R}^p\times
{\mathbb
R}^{p\times m}\rightarrow {\mathbb R}^p\\
&\quad\mbox{being}~{\cal P}\otimes{\cal B}_p\otimes{\cal B}_{p\times m}/{\cal B}_{p}~\mbox{measurable},\\
 g=&(g_{kl})_{p\times d}: \Gamma\times\Omega\times[0,T]\times{\mathbb R}^p\times{\mathbb
R}^{p\times m}\rightarrow {\mathbb
R}^{p\times d}\\
&\quad\mbox{being}~{\cal P}\otimes{\cal B}_p\otimes{\cal B}_{p\times m}/{\cal B}_{p\times d}~\mbox{measurable},\\
h=&(h_{kl})_{p\times m}: \Gamma\times\Omega\times[0,T]\times{\mathbb R}^p\times{\mathbb
R}^{p\times m}\rightarrow {\mathbb
R}^{p\times m}\\
&\quad\mbox{being}~{\cal P}\otimes{\cal B}_p\otimes{\cal B}_{p\times m}/{\cal B}_{p\times m}~\mbox{measurable}.
\end{array}
$$
Moreover, the following conditions are satisfied
$$f(\cdot, 0)\in {\mathbb M}^2(0,T;{\mathbb R}^p), ~ g(\cdot, 0)\in {\mathbb M}^2(0,T; {\mathbb R}^{p\times d})~\mbox{and}~ h(\cdot, 0)\in {\mathbb M}^2(0,T; {\mathbb R}^{p\times m}),$$
and there exists $c>0$ such that
$$
\begin{array}{ll}
|f(t, x_1, y_1)-f(t,x_2, y_2)|&\vee |g(t,x_1, y_1)-g(t,x_2,y_2)|\vee |h(t,x_1, y_1)-h(t,x_2,y_2)|\\
&\leq c(|x_1-x_2|+|y_1-y_2|),\\
\end{array}
\eqno(4.2)
$$
for all $x_1,x_2\in{\mathbb R}^p, y_1,y_2\in{\mathbb R}^{p\times m}, (\gamma,\omega,t)$-a.e.; and there exists $\alpha>0$ such that
$$|h(t, x, y_1)-h(t,x, y_2)|\geq \alpha |y_1-y_2|\eqno(4.3)$$
for all $x\in{\mathbb R}^p, y_1,y_2\in{\mathbb R}^{p\times m}, (\gamma,\omega,t)$-a.e.
Note that (4.3) being satisfied  is the case of  UBSDE  (3.10), with $\alpha=1$. Thus, (4.2) and (4.3) imply that for all $x\in {\mathbb R}^p$ and $(\gamma,\omega,t)$-a.e., the mapping $y\rightarrow h(t,x,y)$ is a bijection from ${\mathbb R}^{p\times m}$ onto itself. In fact, one-to-one property follows at once from (4.2) and the onto property from continuity (4.3), injectivity, and the fact that $\lim_{|y|\rightarrow +\infty}|h(t,x,y)|=+\infty$ for all $x\in{\mathbb R}^p$, $(\gamma,\omega,t)$-a.e. Thus we have the following theorem.

\vskip12pt

{\noindent\bf Theorem 4.1.} {\sl In UBSDE (4.1), if the above conditions on $\xi, f,g,h$ are satisfied, in particular (4.2) and (4.3) hold, then there exists a unique pair $(X, Y)\in {\mathbb M}^2(0,T;{\mathbb R}^p)\times {\mathbb M}^2(0,T;{\mathbb R}^{p\times m})$ which solves  UBSDE (4.1).
}

\vskip12pt

{\noindent\it Proof.} The proof is an adaptation, with essentially obvious changes, of the proofs the previous results.
We just indicate the three steps, and explain the one new argument which is needed in the first step. The first step consists in studying  UBSDE
$$
\begin{array}{ll}
& X(t)+\int_t^Tf(s)ds+\int_t^Tg(s)dC(s)+\int_t^Th(s,Y(s))dB(s)=\xi,
\end{array}
\eqno(4.4)
$$
where $h$ satisfy the simplified version of (4.2) and (4.3) obtained by suppressing the dependence in $x$. The second step solves  UBSDE
$$
\begin{array}{ll}
& X(t)+\int_t^Tf(s, Y(s))ds+\int_t^Tg(s,Y(s))dC(s)+\int_t^Th(s,Y(s))dB(s)=\xi,
\end{array}
$$
where $f,g,h$ satisfies the simplified version of (4.2) and (4.3) obtained by suppressing the dependence in $x$. The third step solves  UBSDE (4.1).

Let us only discuss  UBSDE  (4.4). From Proposition 3.1, there exists a unique pair $(X,\bar Y)$ such that
$$
 X(t)+\int_t^Tf(s)ds+\int_t^Tg(s)dC(s)+\int_t^T\bar{Y}(s)dB(s)=\xi.
$$
It remains only to show that given $\bar{Y}\in {\mathbb M}^2(0,T;{\mathbb R}^{p\times m})$, there exists a unique $Y\in {\mathbb M}^2(0,T;{\mathbb R}^{p\times m})$ such that $h(t,Y(t))=\bar{Y}(t)~ (\gamma,\omega)$-a.e. From the properties of $h$, it follows that for any $(\gamma,\omega, t,y)\in \Gamma\times \Omega\times[0,T]\times {\mathbb R}^{p\times m}$, there exists a unique element $\phi_t(\gamma,\omega,y)$ of ${\mathbb R}^{p\times m}$ such that $h(\gamma,\omega, t,\phi_t(\gamma,\omega, y))=y$. It remains only to show that $\phi$ is ${\cal P}\otimes {\cal B}_{p\times m}/{\cal B}_{p\times m}$ measurable.

Without loss of generality, we can assume that $\Gamma=C([0,T]; {\mathbb R}^d), C_t(\gamma)=\gamma(t)$ and $\Omega=C([0,T]; {\mathbb R}^m)$, $ B_t(\omega)=\omega(t)$. ${\cal L}$ and ${\cal F}$ are the Borel fields over $\Gamma$ and $\Omega$, respectively. Note that the mapping
$$G(\gamma,\omega, t, y)=(\gamma,\omega,t, h(\gamma,\omega,t,y))$$
is a bijection from ${\cal E}=\Gamma\times\Omega\times[0,T]\times{\mathbb R}^{p\times m}$ onto itself. Since ${\cal E}$ is a complete and separable metric space, it follows that from Theorem 10.5, page 506 in Ethier and Kurtz \cite{EK} that $G^{-1}$ is Borel measurable, i.e. ${\cal L}\otimes{\cal F}\otimes{\cal B}[0,T]\otimes{\cal B}_{p\times m}$ measurable. By considering for each $t$ the  restriction of the same map to $C([0,t]; {\mathbb R}^{d})\times C([0,t]; {\mathbb R}^{m})\times [0,t]\times {\mathbb R}^{p\times m}$, we obtain that $G^{-1}$ is ${\cal P}\otimes {\cal B}_{p\times m}$ measurable, which proves the above claim for $\phi$. Thus the proof is complete.\hfill$\Box$

\vskip12pt

\section{Conclusions}
The theory of  UBSDEs combines the theory of backward stochastic differential equations with uncertainty theory, which is a new study field and can be applied in such fields as control theory, physics, engineering, biology, economics and other. It will be an
important tool to deal with uncertain stochastic systems with final values satisfying certain conditions. By preparing three propositions, we have proven an existence and uniqueness theorem of solution to the general UBSDEs under the uniform Lipschitzian
condition, which is the main contribution of this paper. It is believed that our results will be helpful to the study of  UBSDEs.

 In the study of control systems, for example, Shen et al. \cite{SWL} discuss the sampled-data synchronization control problem for a class of dynamical networks. The sampling period is assumed to be time-varying that switches between two different values in a random
way with given probability. The addressed synchronization control problem is first formulated as an exponentially mean-square stabilization problem for a new class of dynamical networks that
involve both the multiple probabilistic interval delays and the sector-bounded nonlinearities. Then, a novel
Lyapunov functional is constructed to obtain sufficient conditions,
under which the dynamical network is exponentially mean-square
stable. In Wang et al. \cite{WSSW}, a class of nonlinear stochastic time-delay network-based systems with probabilistic data missing is investigated. A nonlinear stochastic system with state delays is employed
to model the networked control systems where the measured output and the input signals are quantized by two logarithmic quantizers, respectively. Moreover, the data missing phenomena
are modeled by introducing a diagonal matrix composed of
Bernoulli distributed stochastic variables taking values of 1
and 0, which describes that the data from different sensors
may be lost with different missing probabilities. However, in reality, due to data missing, the parameters estimation is often imprecise so that the system show uncertainty beyond probabilistic uncertainty. So considering only the randomness of the system is inadequate, we think it is necessary to include Liu's uncertainty in the random system. Thus, the theory of uncertain random systems can be applied to the above complex control systems, moreover UBSDEs can be utilized to explore the related nonlinear systems, which is one of our future research direction. On the other hand, we will further explore an existence and uniqueness of solutions to  UBSDEs under non-Lipschitzian condition.

% conference papers do not normally have an appendix

% use section* for acknowledgement

\section*{Acknowledgment} The author is highly grateful to two referees for their valuable comments and suggestions.

% trigger a \newpage just before the given reference
% number - used to balance the columns on the last page
% adjust value as needed - may need to be readjusted if
% the document is modified later
%\IEEEtriggeratref{8}
% The "triggered" command can be changed if desired:
%\IEEEtriggercmd{\enlargethispage{-5in}}

% references section

% can use a bibliography generated by BibTeX as a .bbl file
% BibTeX documentation can be easily obtained at:
% http://www.ctan.org/tex-archive/biblio/bibtex/contrib/doc/
% The IEEEtran BibTeX style support page is at:
% http://www.michaelshell.org/tex/ieeetran/bibtex/
%\bibliographystyle{IEEEtran}
% argument is your BibTeX string definitions and bibliography database(s)
%\bibliography{IEEEabrv,../bib/paper}
%
% <OR> manually copy in the resultant .bbl file
% set second argument of \begin to the number of references
% (used to reserve space for the reference number labels box)
%\begin{thebibliography}{1}

\vspace{1cm}\titleformat{\section}{\large\bfseries}{\S\arabic{section}}{1em}{}

\vskip 24pt
School of Mathematics and Physics, Anhui Polytechnic University, Wuhu 241000, China\\
\indent Email: wyfei@ahpu.edu.cn.

% that's all folks
\end{document}